\documentclass[english]{article}
\usepackage[utf8]{inputenc}
\usepackage[T1]{fontenc}
\usepackage{babel}
\usepackage{amsmath}
\usepackage{graphicx}
\usepackage{fancyhdr}
\usepackage[a4paper,top=3cm,bottom=2cm,left=3cm,right=3cm,marginparwidth=1.75cm]{geometry}

\begin{document}

\title{\bfseries Dynamical Systems on Compact Metrizable Groups}
\author{
    Binghui Xiao \\
    \textit{College of Mathematics and Statistics,} \\
    \textit{Chongqing University, Chongqing, 401331, P.R.China} \\
    Email: \texttt{202106021011@stu.cqu.edu.cn}
}
\date{} 
\maketitle
\thispagestyle{fancy}

\vspace{1cm}

\textbf{Abstract:} This paper aims to extend Kenneth R. Berg's findings on the maximal entropy theorem and the ergodicity of measure convolution to the case of surjective homomorphisms. 

The article further explores dynamical systems under surjective homomorphism in detail, especially the variation of entropy. Let \(G_{1}\) and \(G_{2}\) be compact metrizable groups, and suppose that \(G_{2}\) acts freely on \(G_{1}\), the continuous mapping \(T_{1}\) and homomorphism \(T_{2} :G_{2} \to G_{2}\) satisfy \(T_{1}(yx)=T_{2}(y)T_{1}(x)\), where \(y\in G_{2}, x\in G_{1}\). If \(\mu_{0} \in M(T_{0})\), \(\mu_{0}'\) is the Haar extension of \(\mu_{0}\), we proved that when \(\mu \in M(T_{1},\mu_{0})\), the entropy \(h(T_{1},\mu_{0}')\) is always greater than or equal to \(h(T_{1},\mu)\); if \(\mu_{0}'\) is ergodic with respect to \(T_{1}\), and the Haar measure \(m\) on \(G_{2}\) is ergodic with respect to \(T_{2}\), and if \(h(T_{1},\mu_{0}')<\infty\), then the entropy \(h(T_{1},\mu_{0}')\) is greater than \(h(T_{1},\mu)\).

Finally, this paper also specifically discusses the ergodicity of the convolution of invariant measures. Let \(T\) be a surjective homomorphism on \(G\), if \((G,T,\mathcal{F},\mu)\) and \((G,T,\mathcal{F},\nu)\) are disjoint ergodic dynamical systems, then \(\mu*\nu\) is ergodic. Via a proof by contradiction, the study demonstrates that the measure convolution of two disjoint ergodic dynamical systems can maintain ergodicity under the condition that \(T\) is a surjective homomorphism on \(G\).

\vspace{1cm}
\begin{center}
{\large \bf 1. Introduction}
\end{center}

 In his seminal paper "Convolution of Invariant Measures, Maximal Entropy," Kenneth R. Berg initially outlines two prevalent methods for defining the convolution of Borel measures $\mu$ and $\nu$ on a locally compact Hausdorff group $G$. The first definition is grounded in the Riesz Representation Theorem, considering the unique regular Borel measure that matches specified integral conditions with continuous functions on $G$. The second definition is more direct, explicitly defining convolution for every Borel subset $E$ of $G$. Intriguingly, when $\mu$ and $\nu$ are both regular measures, these definitions are shown to be equivalent, laying a solid foundation for the core discussions that follow in the paper.
 
 Berg's research focuses on compact metrizable groups, delving into the impact of invariant measure convolution on system entropy. Through a series of precise theorems and corollaries, Berg elucidates how the convolution of invariant measures affects the entropy of a system and reveals that under certain conditions, the convolution of two ergodic measures remains ergodic. These insights are of paramount importance for understanding the long-term behavior of dynamical systems, especially when considering the interaction of multiple independent systems.

Furthermore, Berg underscores the unique position of Haar measure in maximizing entropy. He demonstrates that system entropy reaches its maximum with the Haar measure, and if the associated transformation is ergodic and the entropy $h(m)$ is finite, then Haar measure is uniquely characterized as the measure that maximizes entropy. This discovery not only highlights the special property of Haar measure in capturing the system's maximal uncertainty and complexity but also offers a new perspective for understanding and applying dynamical systems.

Building on Berg's contributions, this study aims to further extend his findings, generalizing the main results from isomorphic cases to surjective homomorphism,the main results are theorem 2.2,theorem 3.1,theorem 3.3,theorem 4.1. Our goal is to deepen the understanding of the convolution of invariant measures and entropy maximization problems on compact metrizable groups and explore how these theories can address real-world problems, paving new paths for the theory and application of dynamical systems. By integrating measure theory, ergodic theory, and entropy theory, we hope to provide deeper insights and practical tools for this field, thereby fostering the development and application of dynamical systems theory.

\vspace{0.5cm}{\large\bf notation}

\bigskip 
Let $G$ be a compact metrizable group, $T$ is a surjective homomorphism on $G$, $m$ represents the Haar measure, Let $\mathcal{B}$ be the Borel field. Let $M(T)$ denotes the set of all Borel probability measures on $G$ that are invariant under the map $T$, and $\pi$ represents the projection mapping. Define $\pi : G \times G \to G, \pi(x, y) = xy$, and $\pi_1(x, y) = x, \pi_2(x, y) = y$. We set $\mathcal{B}_c = \pi^{-1}(\mathcal{B}), \mathcal{B}_1 = \pi_1^{-1}(\mathcal{B}), \mathcal{B}_2 = \pi_2^{-1}(\mathcal{B})$. Where $\mathcal{B}$ is the Borel $\sigma$-algebra of $G$.Given fields $\mathcal{A}$, $\mathcal{B}$, $\bigvee_{i \in I} \phi_i$, etc., as in reference [23]. The entropy functions $H(\mathcal{A})$, $H(\mathcal{A}|\mathcal{B})$, and $h(T)$ will also be understood as in [23]. We shall use the notation of a measurable partition, and the correspondence between partitions and fields. A convenient reference is in [23].

\vspace{1cm}
\begin{center}
{\large\bf 2. Entropy of Convolution of Invariant Measures}
\end{center}

\textbf{Lemma 2.1. }${\rm {\mathcal B}}_{i} \vee {\rm {\mathcal B}}_{c} ={\rm {\mathcal B}}\times {\rm {\mathcal B}}$ holds for $i=1,2$.

\textbf{Proof of Lemma 2.1:} We only need to prove that this theorem holds for $i=1$, and the case for $i=2$ can be similarly proved.On one hand, we know that $\pi _{1} $ and $\pi $ are both measurable mappings from $(G\times G,{\rm {\mathcal B}}_{1} \vee {\rm {\mathcal B}}_{c} )$ to $(G,{\rm {\mathcal B})}$. Therefore, $\pi _{1} \times \pi (x,y)=(xy,x)$ is a measurable mapping from $(G\times G,{\rm {\mathcal B}}_{1} \vee {\rm {\mathcal B}}_{c} )$ to $(G\times G,{\rm {\mathcal B}}\times {\rm {\mathcal B}})$.Next, define $\sigma (x,y)=y^{-1} x,$ then $\sigma $ is a measurable mapping from $(G\times G,{\rm {\mathcal B}}\times {\rm {\mathcal B}})$ to $(G,{\rm {\mathcal B})}$, and\[\sigma (\pi _{1} \times \pi (x,y))=\sigma (xy,x)=x^{-1} xy=y=\pi _{2} (x,y).\] This shows that $\pi _{2} $ is a measurable mapping from $(G\times G,{\rm {\mathcal B}}_{1} \vee {\rm {\mathcal B}}_{c} )$ to $(G,{\rm {\mathcal B})}$.On the other hand, by the definition of $\pi _{2} $, it is known that $\pi _{2} $ is a measurable mapping from $(G\times G,{\rm {\mathcal B}}\vee {\rm {\mathcal B}})$ to $(G,{\rm {\mathcal B})}$, therefore ${\rm {\mathcal B}}\times {\rm {\mathcal B}}={\rm {\mathcal B}}_{1} \vee {\rm {\mathcal B}}_{c} .$

Lemma 2.1 is thus proved.     

\bigskip 
\textbf{Lemma 2.2. }For probability measures $\mu $ and $\nu $ on group G, if $\mu ,\nu $ are invariant measures with respect to group homomorphisms $T,S$ respectively, then $h_{\mu \times \nu } (T\times S)=h_{\mu } (T)+h_{\nu } (S)$.

\textbf{Proof of Lemma 2.2:} Let $\alpha =\{ A_{i} \} _{i=1}^{n} ,\beta =\{ B_{j} \} _{j=1}^{m} $ be measurable partitions of $G$, then
\begin{align*}
H_{\mu \times \nu } (T\times S) &= -\sum _{i=1}^{n}\sum _{j=1}^{m}\mu \times \nu (A_{i} \times B_{j} )\log \mu \times \nu (A_{i} \times B_{j} ) \\
&= -\sum _{i=1}^{n}\sum _{j=1}^{m}\mu (A_{i})\nu (B_{j})\log \mu (A_{i})\nu (B_{j}) \\
&= -\sum _{i=1}^{n}\sum _{j=1}^{m}\mu (A_{i})\nu (B_{j})\log \mu (A_{i}) -\sum _{i=1}^{n}\sum _{j=1}^{m}\mu (A_{i})\nu (B_{j})\log \nu (B_{j}) \\
&= -\sum _{i=1}^{n}\mu (A_{i})\log \mu (A_{i}) -\sum _{j=1}^{m}\nu (B_{j})\log \nu (B_{j}) \\
&= H_{\mu } (T)+H_{\nu } (S)
\end{align*}

Following the definition of measure-theoretic entropy \[h_{\mu } (T)={\mathop{\sup }\limits_{\alpha \in {\rm {\mathcal I}}}} {\mathop{\lim }\limits_{n\to \infty }} \frac{1}{n} H_{\mu } (\mathop{{\rm \vee }}\limits_{i=0}^{n-1} T^{-i} \alpha ).\]

Lemma 2.2 is thus proved.

\bigskip
\textbf{Lemma 2.3.} For a group homomorphism $T: G \to G$ on a compact metrizable group $G$ and probability measures $\mu, \nu \in M(T)$, we have
\[h_{\mu * \nu}(T) \leq h_{\mu}(T) + h_{\nu}(T)\].

\textbf{Proof of Lemma 2.3: }Since $(T\times T,{\rm {\mathcal B}}_{c} ,\mu \times \nu )$ and $(T,{\rm {\mathcal B}},\mu *\nu )$ are two dynamical systems

\noindent conjugate to each other, this can be demonstrated by explicitly constructing a homeomorphism $\pi :G\times G\to G,$ such that $\pi \circ (T\times T)=T\circ \pi ,$ since the measure entropy of dynamical systems is invariant under conjugacy, we have 
\[h_{\mu *\nu } (T,{\rm {\mathcal B}})=h_{\mu \times \nu } (T\times T,{\rm {\mathcal B}}_{c} ).\] 

On the other hand, from the lemma, since ${\rm {\mathcal B}}_{i} \vee {\rm {\mathcal B}}_{c} ={\rm {\mathcal B}}\times {\rm {\mathcal B}}$, we have
\[h_{\mu \times \nu } (T\times T,{\rm {\mathcal B}}_{c} )\le h_{\mu \times \nu } (T\times T,{\rm {\mathcal B}}_{c} \vee {\rm {\mathcal B}}_{i} )=h_{\mu \times \nu } (T\times T,{\rm {\mathcal B}}\times {\rm {\mathcal B}}).\] 

Finally, from Lemma 2.2, we obtain
\[h_{\mu \times \nu } (T\times T,{\rm {\mathcal B}}\times {\rm {\mathcal B}})=h_{\mu } (T,{\rm {\mathcal B}})+h_{\nu } (T,{\rm {\mathcal B}}).\] 

Combining the above points, we derive the inequality $h_{\mu *\nu } (T)\le h_{\mu } (T)+h_{\nu } (T).$

Lemma 2.3 is thus proved.                          

\bigskip
\textbf{Corollary 2.1.} For a group homomorphism $T: G \to G$ on a compact metric group $G$ and a probability measure $\mu \in M(T)$, we have
\[h_m(T) \ge h_{\mu}(T).\]

\textbf{Proof of Corollary 2.1:} Since \begin{align*}
m*\mu (E) &= m \times \mu (\pi ^{-1} (E))\\
&= m \times \mu (\{(x,y) \,|\, (x,y) \in \pi ^{-1} (E)\})\\
&= m \times \mu (\{(x,y) \,|\, xy \in E\})\\
&= \int_{G} \int_{Ey^{-1}} dm(x) d\mu(y)\\
&= m(Ey^{-1})\\
&= m(E)
\end{align*}

Therefore, we have $m*\mu =m$. Assume $h_{m} (T)<\infty $, thus $h_{\mu } (T)\le h_{m*\mu } (T)=h_{m} (T).$

Corollary 2.1 is thus proved.  

\bigskip
\textbf{Theorem 2.1.} For a group homomorphism $T: G \to G$ on a compact metrizable group $G$ and probability measures $\mu, \nu \in M(T)$, if $h_{\mu}(T) < \infty$, then $h_{\nu}(T) \le h_{\mu * \nu}(T)$; if $h_{\nu}(T) < \infty$, then $h_{\mu}(T) \le h_{\mu * \nu}(T).$

\textbf{Proof of Theorem 2.1:} Let $h_{\mu \times \nu } (T\times T,{\rm {\mathcal B}}\times {\rm {\mathcal B}})$ denote the entropy of the measure $\mu \times \nu $ under the transformation $T\times T$.

First, since $G$ is a compact metrizable group, hence $G$ is separable, thus ${\rm {\mathcal B}}$ is also separable. Therefore, there exists an increasing sequence of finite algebras ${\rm {\mathcal B}}_{n} $ and ${\rm {\mathcal C}}_{n} $, such that $\bigvee_{n=1}^{\infty } {\rm {\mathcal B}}_{n} ={\rm {\mathcal B}}_{c} $ and $\bigvee_{n=1}^{\infty } {\rm {\mathcal C}}_{n} ={\rm {\mathcal B}}_{1} $. For an algebra ${\rm {\mathcal D}}$, we define ${\rm {\mathcal D}}^{-} =\bigvee_{i=1}^{\infty } T^{-i} {\rm {\mathcal D}}$. Thus, we have
\[H_{\mu \times \nu } ({\rm {\mathcal B}}_{n} \vee {\rm {\mathcal C}}_{n} |({\rm {\mathcal B}}_{n} \vee {\rm {\mathcal C}}_{n} )^{-} )\le H_{\mu \times \nu } ({\rm {\mathcal B}}_{n} |{\rm {\mathcal B}}_{n} {}^{-} )+H_{\mu \times \nu } ({\rm {\mathcal C}}_{n} |{\rm {\mathcal C}}_{n} {}^{-} )\] 
As $n\to \infty $, we obtain
\[h_{\mu \times \nu } (T\times T,{\rm {\mathcal B}}\times {\rm {\mathcal B}})\le h_{\mu \times \nu } (T\times T,{\rm {\mathcal B}}_{c} )+h_{\mu \times \nu } (T\times T,{\rm {\mathcal B}}_{1} )\] 

Based on the equivalence relation, we have $h_{\mu } (T)+h_{\nu } (T)\le h_{\mu *\nu } (T)+h_{\mu } (T)$. Specifically, if $h_{\mu } (T)<\infty $, we obtain $h_{\nu } (T)\le h_{\mu *\nu } (T)$. Similarly, if $h_{\nu } (T)<\infty $, we can derive $h_{\mu } (T)\le h_{\mu *\nu } (T).$

Theorem 2.1 is thus proved. 

\bigskip
\textbf{Corollary 2.2.} If $h_{\mu}(T) = 0$, then $h_{\mu * \nu}(T) = h_{\nu}(T).$

\textbf{Proof of Corollary 2.2:} Since $h_{\mu } (T)=0$, we have $0\le h_{\mu *\nu } (T)-h_{\nu } (T)$. Therefore, $ h_{\nu } (T)= h_{\mu } (T)+ h_{\nu } (T)\ge h_{\mu *\nu } (T)\ge h_{\nu } (T)$.

Corollary 2.2 is thus proved.

\bigskip
\textbf{Theorem 2.2.}With the previous notations,
\begin{enumerate}
\item[(1)] If $T$ is a surjective homomorphism on $G$, then $m \in M(T)$.
\item[(2)] If $T$ is a homomorphism on $G$, and $\mu, \nu \in M(T)$, then $\mu * \nu \in M(T)$.
\end{enumerate}

\bigskip 
\textbf{Proof of Theorem 2.2:} 
(1) Since $T$ is a homomorphism and surjective, for each $f\in C_{c}(G)$, $\forall g\in G$, we have
\begin{align}
\int_{G}f(gx)d(m\circ T^{-1})(x) &= \int_{G}f(g\cdot Tx)dm(x) \tag{2.1} \label{eq:custom1} \\
    &= \int_{G}f\circ T(h\cdot x)dm(x) \tag{2.2} \label{eq:custom2} \\
    &= \int_{G}f\circ T(x)dm(x) \notag \\
    &= \int_{G}f(x)d(m\circ T^{-1})(x) \notag
\end{align}

At (2.1), let $T^{-1}x=y\Rightarrow x=TT^{-1}x=Ty$.

At (2.2), since $T$ is surjective, so $\forall g\in G,\exists h\in G,\mathrm{s.t\;} Th=g$.

Therefore, according to the uniqueness of the Haar measure on compact metric groups, $m\circ T^{-1} =cm$. On the other hand, for the Borel set $E_{0}$, we have $0<m\circ T^{-1}(E_{0})=m(E_{0})<\infty$, thus $c=1$, so $m\in M(T)$.

(2) Since $T$ is a group homomorphism, then\[T(\pi (x,y))=T(xy)=Tx\cdot Ty=\pi (Tx,Ty)=(\pi (T\times T))(x,y)\]
so $T\pi =\pi (T\times T)$. If $\mu ,\nu \in M(T)$, for any $E\in \mathcal{B}$, we have\[\mu *\nu (T^{-1} E)=\mu \times \nu (\pi ^{-1} (T^{-1} E))=\mu \times \nu ((T^{-1} \times T^{-1} )\pi ^{-1} E))=\mu \times \nu (\pi ^{-1} E)=\mu *\nu (E)\]

On the other hand, if $\mu ,\nu \in M(T)$, choose Borel sets $E_{1} ,E_{2}$ such that $0<\mu (E_{1})<\infty$ and $0<\nu (E_{2})<\infty$, and let $E_{0} =\pi (E_{1} \times E_{2})$, then
\[\mu *\nu (E_{0})=\mu \times \nu (\pi ^{-1} E_{0})=\mu \times \nu (\pi ^{-1} \pi (E_{1} \times E_{2}))\ge \mu \times \nu (E_{1} \times E_{2})>0\]

Theorem 2.2 is thus proved.

\vspace{1cm}
\begin{center}
{\large\bf 3. Group extensions}
\end{center}

This chapter delves into the principle of maximum entropy and the convolution ergodicity of invariant measures on compact metrizable groups. This section first introduces some definitions, standardizes symbols, and then presents the main results of this chapter.

Let $\mathcal{I}$ denote the set of measurable partitions of $X$ with finite entropy. Let $\alpha_T = \bigvee_{n=-\infty}^{n=\infty} T^n \alpha$, $\alpha^- = \bigvee_{i=1}^{\infty} T^{-i} \alpha, \alpha_n = \bigvee_{i=1}^{n} T^i \alpha, \alpha^n = \bigvee_{i=0}^{n} T^{-i} \alpha$, $\mathcal{A}$ is the algebra corresponding to the partition $\alpha_T$, $h(T,\alpha)$ is the entropy of the system $(X, \mathcal{A}, T, \mu)$. Similarly define $\beta_T, \mathcal{B}$, and $h(T, \beta)$.

Let $G_1$ be a compact metrizable group, and let $G_2$ be a closed normal subgroup of $G_1$, with the action of $G_2$ on $G_1$ being left multiplication, then the quotient group $G_1 / G_2$ is also a compact metrizable group, with $\varphi: G_1 \to G_1 / G_2$ representing the natural projection. The continuous mapping $\pi: G_2 \times G_1 \to G_1$ such that $\pi(y,x) = yx (y \in G_2, x \in G_1)$, $T_1: G_1 \to G_1$ is a homeomorphism and $T_2: G_2 \to G_2$ is an isomorphism, and satisfy $T_1(yx) = T_2(y)T_1(x)$, where $y \in G_2, x \in G_1$. In this case, $T_1$ induces an isomorphic mapping $T_0$ on $G_1 / G_2$, defined as $T_0([x]) = [T_1(x)]$, where $[x]$ represents the equivalence class that element $x$ belongs to.

$T_{0} $ is well-defined. Indeed, if $[x_{1} ]=[x_{2} ]$, there exists $y\in G_{2} $ such that $x_{1} =yx_{1} $. Applying $T_{1} $ to both sides yields$T(x_{1} )=T_{1} (yx_{2} )=T_{2} (y)T_{1} (x_{2} ),{\rm \; \; }T_{2} (y)\in G_{2} $,therefore $[T_{1} (x_{1} )]=[T_{1} (x_{2} )]$.If $\mu _{0} \in M(T_{0} )$, $M(T_{1} ,\mu _{0} )=\{ \mu :\mu \in M(T_{1} ),\mu \circ \varphi ^{-1} =\mu _{0} \} $, $m$ is the normalized Haar measure on $G_{2} $, the Haar extension of $\mu _{0} $, denoted $\mu _{0} '$, is defined by
\[\int _{G_{1} }fd\mu _{0} ' =\int _{G_{1} {\rm / }G_{2} }\int _{G_{2} }f(yx) dm(y)d\mu _{0}  (x)\] 
for $f\in C_{c} (G)$.

\noindent Because
\begin{align}
\int _{G_{1} / G_{2} }fd\mu _{0} ' \circ \varphi ^{-1} &= \int _{G_{1} / G_{2} }f\circ \varphi d \mu _{0} ' \notag\\
&= \int _{G_{1} / G_{2} }\int _{G_{2} }f(\varphi (yx))dm(y)d\mu _{0} (x) \notag\\
&= \int _{G_{1} / G_{2} }\int _{G_{2} }f([x])dm(y)d\mu _{0} (x) \notag\\
&= \int _{G_{1} / G_{2} }f([x])d\mu _{0} (x)\notag
\end{align}

So $\mu _{0} '\circ \varphi ^{-1} =\mu _{0} $, therefore $\mu _{0} '\in M(T_{1} ,\mu _{0} )$.

\bigskip
\textbf{Definition 3.1.} If $G_{1}$ and $G_{2}$ is a compact metrizable group, and $T_{1} :G_{1} \to G_{1}$, $T_{2} :G_{2} \to G_{2}$ are continuous mappings, if $G_{2}$ is a compact metrizable group acting on $G_{1}$, if $T_{1}(yx)=T_{2}(y)T_{1}(x)$, where $y \in G_2, x \in G_1$,then we say $T_{1}$ $G_{2}{\rm-}$ commutes with $T_{1}$.

\bigskip
\textbf{Definition 3.1.} The the group $G_{2}$ act freely on $G_{1}$, if for any element $x$ in $G_{1}$ and any element $g$ in $G_{2}$ other than the identity element $e$, we have $g\cdot x\ne x$.

\bigskip
\textbf{Lemma 3.1.} Assume $\mu_{0} \in M(T_{0})$ is ergodic, then
\begin{enumerate}
\item[(1)] $M(T_{1}, \mu_{0})$ is a convex set;

\item[(2)] $M(T_{1}, \mu_{0})$ is closed under the weak* topology;

\item[(3)] The extreme points in $M(T_{1}, \mu_{0})$ are the ergodic measures with respect to $T_{1}$.
\end{enumerate}

\textbf{Proof of Lemma 3.1:}(1) Let $\mu, \nu \in M(T_{1}, \mu_{0}), a+b=1$, for $\forall B \in {\mathcal B}$, because
\begin{align*}
(a\mu + b\nu) \circ \varphi^{-1}(B) &= a\mu \circ \varphi^{-1}(B) + b\nu \circ \varphi^{-1}(B) \\
&= a\mu_{0}(B) + b\mu_{0}(B) \\
&= \mu_{0}(B)
\end{align*}

Therefore, $a\mu + b\nu \in M(T_{1}, \mu_{0})$, hence $M(T_{1}, \mu_{0})$ is a convex set.

\bigskip
(2) The definition of weak* convergence in measure spaces is that for every $f\in C_{c} (G)$, we have

$$\int _{G}fd\mu _{n} \to \int _{G}fd\mu \, (n\to \infty ).$$

In this case, we have a sequence of measures $\mu _{n} \in M(T_{1} ,\mu _{0} )$, and $\mu _{n} \to \mu$, each $\mu _{n}$ satisfies $\mu _{n} \circ \varphi ^{-1} =\mu _{0}$.

Now, we want to prove that $\mu \circ \varphi ^{-1} =\mu _{0}$. To prove this, we can consider for $\forall f\in C_{c} (G_{1})$, since $\mu _{n} \circ \varphi ^{-1} =\mu _{0}$, we have
\[\int _{G_{1} }f (\varphi (x))d\mu _{n} (x)=\int _{G_{1} }f (x)d\mu _{0} (x)\]

for all $n$.

Since $\mu _{n}$ weakly* converges to $\mu$, we have
\[\int _{G_{1} }f (\varphi (x))d\mu _{n} (x)\to \int _{G_{1} }f (\varphi (x))d\mu (x)\]

as $n\to \infty$.

Combining these two equations, we get
\[\int _{G_{1} }f (x)d\mu _{0} (x)=\int _{G_{1} }f (\varphi (x))d\mu (x)\]

for all $f\in C_{c} (G_{1})$.

This means that $\mu \circ \varphi ^{-1} =\mu _{0}$, which completes the proof.

(3) We know that the extreme points of $M(T_{1})$ are precisely the ergodic measures in $M(T_{1})$. Thus, an ergodic measure in $M(T_{1}, \mu _{0})$ is an extreme point of $M(T_{1})$, and therefore also an extreme point of $M(T_{1}, \mu _{0})$.Conversely, suppose $\mu$ is an extreme point of $M(T_{1}, \mu_{0})$. We want to prove that $\mu$ is also an extreme point of $M(T_{1})$. Suppose $\mu = a\tau + b\nu$, where $\tau, \nu \in M(T_{1})\; (a, b > 0, a + b = 1)$. Then$\mu_{0} = \mu \circ \varphi^{-1} = a\tau \circ \varphi^{-1} + b\nu \circ \varphi^{-1}$.Since $\mu_{0}$ is ergodic, we have $\tau \circ \varphi^{-1} = \nu \circ \varphi^{-1} = \mu_{0}$. Therefore, $\tau, \nu \in M(T_{1}, \mu_{0})$, and so $\tau = \nu = \mu$. Hence, $\mu$ is an extreme point of $M(T_{1})$.

Lemma 3.1 is thus proved.

\bigskip
\textbf{Corollary 3.1.} If $\mu_{0} \in M(T_{0})$ is ergodic, and for all ergodic measures $\mu$ in $M(T_{1}, \mu_{0})\; (\mu \neq \mu_{0}')$, we have
\[h(T_{1}, \mu) < h(T_{1}, \mu_{0}')\]
Then this inequality also holds for all $\mu \in M(T_{1}, \mu_{0})\; (\mu \neq \mu_{0}')$, and $\mu_{0}'$ itself is ergodic.

\textbf{Proof of Corollary 3.1:} Let $\mu \in M(T_{1}, \mu_{0})\; (\mu \neq \mu_{0}')$. According to Choquet's theorem, $\mu$ can be represented as $\mu = \int \nu d\lambda$, where $\lambda$ is a probability measure on $M(T_{1}, \mu_{0})$ with support on the extreme points. Then
\[h(T_{1}, \mu) = \int h(T_{1}, \nu) d\lambda\]
and by using Corollary 4.2.1, if $\mu \neq \mu_{0}'$, then $h(T_{1}, \mu) < h(T_{1}, \mu_{0}')$. By similar reasoning, if $\mu_{0}' = \int \nu d\lambda$, then $\nu = \mu_{0}'$ almost everywhere, and $\mu_{0}'$ is ergodic.

Corollary 3.1 is thus proved.

\bigskip
Next is a formula for calculating the entropy $h(T_{1}, \mu_{0}')$

\textbf{Lemma 3.2.${}^{\ }$${}^{[}$${}^{18}$${}^{]}$ (The Entropy Addition Theorem)} Let $(M, \mathcal{B}, \mu)$ be a Lebesgue space, which is a product of a Lebesgue space $(X, \mathcal{F}, \nu)$ and a compact, separable group $G$ with Borel sets and Haar measure $m$. All measures are normalized, i.e., $\mu(M) = \nu(X) = m(G) = 1$. Let $T$ be a measure-preserving transformation defined by
\[T(x,g) = (Sx, \sigma(g)\varphi(x))\]
where $S$ is a measure-preserving transformation on $X$, $\sigma$ is a group automorphism on $G$, and $\varphi: X \to G$ is a measurable mapping. In this case, we have the following entropy addition theorem:
\[h(T) = h(S) + h(\sigma)\]

The implication of this theorem is that if a measure-preserving transformation can be decomposed into two parts, one part being a measure-preserving transformation and the other part being an automorphism of a topological group, then the entropy of this transformation can be obtained by calculating the entropies of these two parts and adding them together.

\bigskip
\textbf{Lemma 3.3.} When $T_{1}, T_{2}$ satisfy $T_{1}(yx)=T_{2}(y)T_{1}(x)$, and $G_{2}$ acts freely on $G_{1}$, that is, if $y \neq e$ then $yx \neq x\; (y \in G_{2}, x \in G_{1})$. We have
$$h(T_{1}, \mu_{0}') = h(T_{0}, \mu_{0}) + h(T_{2}, m)$$

\noindent where $\mu_{0} \in M(T_{0})$ and $m$ denotes the Haar measure on $G_{2}$.
Define the mapping $\phi: G_{1} / G_{2} \times G_{2} \to G_{1} / G_{2} \times G_{2}$ such that $\phi(x_{0}, y) = (T_{0}(x_{0}), T_{2}(y)\eta(x_{0}))$, where $\eta: G_{1} / G_{2} \to G_{2}$ is a Borel measurable mapping. By Lemma 3.2, we have:
\[h(\phi, \mu_{0}') = h(T_{0}, \mu_{0}) + h(T_{2}, m)\]

Since $T_{1}$ and the mapping $\phi$ are Borel isomorphic, i.e., there exists a Borel measurable bijection $\psi: G_{1} \to G_{1} / G_{2} \times G_{2}, \psi(x) = (\varphi(x), \eta(\varphi(x)))$, such that the following equation holds:
\[\psi \circ T_{1} = (T_{0}(x_{0}), T_{2}(y)\phi(x_{0})) \circ \psi\]
Therefore, $h(T_{1}, \mu_{0}') = h(T_{1}, \phi)$, thus $h(T_{1}, \mu_{0}') = h(T_{0}, \mu_{0}) + h(T_{2}, m)$ holds.

\bigskip
\textbf{Lemma 3.4.} Let $\mu_{0} \in M(T_{0})$, then
\begin{enumerate}
\item[(1)] If $\nu \in M(T_{2}), \mu \in M(T_{1}, \mu_{0})$, then $\nu * \mu \in M(T_{1}, \mu_{0}),$

\item[(2)] If $\mu \in M(T_{1}, \mu_{0})$, then $m * \mu = \mu_{0}'.$
\end{enumerate}

\textbf{Proof of Lemma 3.4:}(1) Since that
\begin{align*}
(\pi \circ (T_{2} \times T_{1}))(y,x) &= \pi(T_{2} y, T_{1} x) \\
&= T_{2}(y)T_{1}(x) \\
&= T_{1}(yx) \\
&= T_{1} \circ \pi(y,x).
\end{align*}

Therefore, $\pi \circ (T_{2} \times T_{1}) = T_{1} \circ \pi$. Therefore, $\nu \in M(T_{2}), \mu \in M(T_{1}) \Rightarrow \nu * \mu \in M(T_{1})$, since $\mu \in M(T_{1}, \mu_{0})$, so for any set $B$ on $G_{1} / G_{2}$, we have $\mu_{0}(B) = \mu(\varphi^{-1} B)$.

\begin{align*}
    (\nu *\mu )\circ \varphi^{-1} (B) &= \nu \times \mu (\pi^{-1} (\varphi^{-1} B)) \\
    &= \nu \times \mu (\{(y,x) \,|\, yx\in \varphi^{-1} B\}) \\
    &= \mu (y^{-1} (\varphi^{-1} B))
\end{align*}

Considering the set $\{ y^{-1} (\varphi^{-1} B)\}$ in $G_{1}$, since $y^{-1}$ is the inverse operation of elements in $G_{2}$, it does not change the equivalence class of elements in $\varphi^{-1} B$, thus $\mu (y^{-1} (\varphi^{-1} B))=\mu_{0} (B)$, hence we obtain $(\nu *\mu )\circ \varphi^{-1} =\mu_{0}$, i.e., $\nu *\mu \in M(T_{1}, \mu_{0})$, completing the proof.

\noindent (2) By definition, it is obvious.When $G_{1} =G_{2}, T_{1} =T_{2}$, then $M(T_{1}, \mu_{0})=M(T), \mu_{0}'=m$, Lemma 1.1 in reference [1] is a special case of this Lemma 3.4.Next, we define the mappings
\[\begin{array}{l} {\pi :G_{2} \times G_{1} \to G_{2}, \pi (x,y)=xy} \\ {\pi_{1} :G_{2} \times G_{1} \to G_{2}, \pi_{1} (x,y)=x} \\ {\pi_{2} :G_{2} \times G_{1} \to G_{1}, \pi_{2} (x,y)=y} \end{array}\] 

Let ${\rm {\mathcal B}}$ be the $\sigma$-algebra generated by $G_{1}$, ${\rm {\mathcal F}}$ be the $\sigma$-algebra generated by $G_{2}$, and let ${\rm {\mathcal B}}_{1} =\pi_{1}^{-1} ({\rm {\mathcal F}})$ and ${\rm {\mathcal B}}_{c} =\pi^{-1} ({\rm {\mathcal F}})$, then ${\rm {\mathcal B}}_{1}$ and ${\rm {\mathcal B}}_{c}$ are sub-$\sigma$-algebras of ${\rm {\mathcal F}}\times {\rm {\mathcal B}}$.

Lemma 3.4 is thus proved.

\bigskip
\textbf{Lemma 3.5.} ${\rm {\mathcal B}}_{1} {\rm \vee {\mathcal B}}_{c} $=${\rm {\mathcal F}}\times {\rm {\mathcal B}}$

\textbf{Proof of Lemma 3.5:} Similar to the proof of Lemma 2.1, on one hand, since $\pi$ is a measurable mapping from $(G_{2} \times G_{1}, {\rm {\mathcal B}}_{1} \vee {\rm {\mathcal B}}_{c})$ to $(G_{1}, {\rm {\mathcal B}})$, and $\pi_{1}$ is a measurable mapping from $(G_{2} \times G_{1}, {\rm {\mathcal B}}_{1} \vee {\rm {\mathcal B}}_{c})$ to $(G_{2}, {\rm {\mathcal F}})$, therefore $\pi \times \pi_{1}(x,y)=(xy,x)$ is a measurable mapping from $(G_{2} \times G_{1}, {\rm {\mathcal B}}_{1} \vee {\rm {\mathcal B}}_{c})$ to $(G_{2} \times G_{1}, {\rm {\mathcal F}}\times {\rm {\mathcal B}})$.Define $\sigma(x,y)=y^{-1}x,$ then $\sigma$ is a measurable mapping from $(G_{2} \times G_{1}, {\rm {\mathcal F}}\times {\rm {\mathcal B}})$ to $(G_{1}, {\rm {\mathcal B}})$, and
\[\sigma(\pi_{1} \times \pi(x,y))=\sigma(xy,x)=x^{-1}xy=y=\pi_{2}(x,y).\]

This shows that $\pi_{2}$ is a measurable mapping from $(G_{2} \times G_{1}, {\rm {\mathcal B}}_{1} \vee {\rm {\mathcal B}}_{c})$ to $(G_{1}, {\rm {\mathcal B}})$.On the other hand, by the definition of $\pi_{2}$, it is known that $\pi_{2}$ is a measurable mapping from $(G_{2} \times G_{1}, {\rm {\mathcal F}}\times {\rm {\mathcal B}})$ to $(G_{1}, {\rm {\mathcal B}})$, therefore ${\rm {\mathcal F}}\times {\rm {\mathcal B}}={\rm {\mathcal B}}_{1} \vee {\rm {\mathcal B}}_{c}.$

Lemma 3.5 is thus proved.

\bigskip
\textbf{Lemma 3.6.} If $\alpha \beta \in {\rm {\mathcal I}}$ and $\alpha \le \beta $, then
\begin{equation}
H(\alpha ^{n} |T^{-(n-1)} \alpha )=\sum _{k=0}^{n-1}H(\alpha |\alpha ^{-}  \vee T^{-1} \beta ^{k} ) \tag{3.1}
\end{equation}

In particular,
\begin{equation}
H(\beta ^{n} |T^{-(n-1)} \beta ^{-} )=nh(T,\beta ).\tag{3.2}
\end{equation}

\textbf{Proof of Lemma 3.6:} Since $\beta ^{k-1} =\beta {\rm \vee }T^{-1} \beta ^{k-1} $, therefore
\[H(\beta ^{k} |T^{-(k-1)} \alpha ^{-} )=H(T^{-1} \beta ^{k-1} |T^{-(k-1)} \alpha ^{-} )+H(\beta |T^{-1} \beta ^{k-1} {\rm \vee }T^{-(k-1)} \alpha ^{-} )\] 
Since $T$ is homomorphic, and $\alpha \le \beta \Rightarrow \alpha {\rm \vee }\beta =\beta $, then
\begin{align*}
T^{-(k-1)} \alpha^{-} \vee T^{-1} \beta^{k-1} &= T^{-(k-1)} \left( \bigvee_{i=1}^{\infty} T^{-i} \alpha \right) \vee T^{-1} \left( \bigvee_{i=0}^{k-1} T^{-i} \beta \right) \\
&= \left( \bigvee_{i=k}^{\infty} T^{-i} \alpha \right) \vee \left( \bigvee_{i=1}^{k} T^{-i} \beta \right) \\
&= \left( \bigvee_{i=k}^{\infty} T^{-i} \alpha \right) \vee \left( \left( \bigvee_{i=1}^{k} T^{-i} \beta \right) \vee \left( \bigvee_{i=1}^{k} T^{-i} \alpha \right) \right) \\
&= \left( \bigvee_{i=1}^{\infty} T^{-i} \alpha \right) \vee \left( \bigvee_{i=1}^{k} T^{-i} \beta \right) \\
&= \alpha^{-} \vee T^{-1} \beta^{k-1}
\end{align*}

Therefore, by induction, we have
\begin{align*}
H(\beta^{k} |T^{-(k-1)} \alpha^{-}) &= H(\beta^{k-1} |T^{-(k-2)} \alpha^{-}) + H(\beta |\alpha^{-} \vee T^{-1} \beta^{k-1}) \\
&= H(\beta |\alpha^{-} \vee T^{-1} \beta^{k-1}) + H(\beta |\alpha^{-} \vee T^{-1} \beta^{k-2}) + \cdots + H(\beta |\alpha^{-} \vee T^{-1} \beta^{0}) \\
&= \sum_{i=0}^{k-1}H(\beta |\alpha^{-} \vee T^{-1} \beta^{i})
\end{align*}

Lemma 3.6 is thus proved.

\bigskip
\textbf{Lemma 3.7.} If $\alpha \beta \in {\rm {\mathcal I}}$, and $\alpha \le \beta ,H(\beta |\alpha^{-} )<\infty $, $T$ is a homomorphism on $G$, then
\begin{equation}
\frac{1}{n} H(\alpha^{n} |T^{-(n-1)} \beta^{-} )\rightarrow h(T,\alpha ) \tag{3.3}
\end{equation}

\textbf{Proof of Lemma 3.7:} Let $\delta$ be a positive number, by Lemma 3.6's (3.2), we get
\[\frac{1}{n} H(\alpha^{n} |T^{-(n-1)} \beta^{-} )\le \frac{1}{n} H(\alpha^{n} |T^{-(n-1)} \alpha^{-} )=h(T,\alpha)\]
Because $\beta^{n} \nearrow T\beta^{-} $, thus
\[T\alpha^{-} \vee \beta^{n} \nearrow T\alpha^{-} \vee T\beta^{-} =T\beta^{-} \Rightarrow \alpha^{-} \vee T^{-1} \beta^{n} \nearrow \beta^{-} \]
Therefore
\[H(\beta |\alpha^{-} \vee T^{-1} \beta^{n} )\searrow H(\beta |\beta^{-} )=h(T,\beta).\]
So, by Lemma 3.6's (3.1), we have $\frac{1}{n} H(\beta^{n} |T^{-(n-1)} \alpha^{-} )=H(\beta |\alpha^{-} \vee T^{-1} \beta^{n} )\searrow h(T,\beta)$

\noindent Next, we only need to prove
\[\frac{1}{n} H(\alpha^{n} |T^{-(n-1)} \beta^{-} )>h(T,\alpha )-\delta \]
holds for all $n$

\noindent Because
\[\frac{1}{n} H(\beta^{n} |T^{-(n-1)} \alpha^{-} )<h(T,\beta )+\delta \]
Thus
\begin{align*}
\frac{1}{n} H(\alpha^{n} |T^{-(n-1)} \beta^{-} ) &= \frac{1}{n} H(\beta^{n} |T^{-(n-1)} \beta^{-} )-\frac{1}{n} H(\beta^{n} |\alpha^{n} \vee T^{-(n-1)} \beta^{-} ) \\
&\ge h(T,\beta )-\frac{1}{n} H(\beta^{n} |\alpha^{n} \vee T^{-(n-1)} \alpha^{-} ) \\
&\ge \frac{1}{n} H(\beta^{n} |T^{-(n-1)} \alpha^{-} )-\delta -\frac{1}{n} H(\beta^{n} |\alpha^{n} \vee T^{-(n-1)} \alpha^{-} ) \\
&= \frac{1}{n} H(\alpha^{n} \vee \beta^{n} |T^{-(n-1)} \alpha^{-} )-\delta -\frac{1}{n} H(\beta^{n} |\alpha^{n} \vee T^{-(n-1)} \alpha^{-} ) \\
&= \frac{1}{n} H(\alpha^{n} |T^{-(n-1)} \alpha^{-} )-\delta \\
&= h(T,\alpha )-\delta.
\end{align*}

Lemma 3.7 is thus proved.

\bigskip
\textbf{Lemma 3.8. }If $\alpha ,\beta \in {\rm {\mathcal I}}$, $T$ is an isomorphism on $G$, then
\begin{equation}
h(\alpha {\rm \vee }\beta )=H(\alpha |\alpha ^{-} {\rm \vee }\beta _{T} )+h(\beta )\tag{3.4}  
\end{equation}

\textbf{Proof of Lemma 3.8:} Because
\begin{align*}
h(T,\beta \vee \alpha) &= \frac{1}{n} H(\beta^n \vee \alpha^n |T^{-(n-1)} \beta^{-} \vee T^{-(n-1)} \alpha^{-}) \\
&= \frac{1}{n} H(\beta^n |T^{-(n-1)} \beta^{-} \vee T^{-(n-1)} \alpha^{-}) + \frac{1}{n} H(\alpha^n |\beta^n \vee T^{-(n-1)} \beta^{-} \vee T^{-(n-1)} \alpha^{-}) \\
&= \frac{1}{n} H(\beta^n |T^{-(n-1)} \beta^{-} \vee T^{-(n-1)} \alpha^{-}) + \frac{1}{n} H(\alpha^n |T\beta^{-} \vee T^{-(n-1)} \alpha^{-})
\end{align*}

\begin{align*}
H(\alpha^k |T\beta^{-} \vee T^{-(k-1)} \alpha^{-}) &= H(T^{-(k-1)} \alpha |T\beta^{-} \vee T^{-(k-1)} \alpha^{-}) \\
&\quad + H(\alpha^{k-1} |T\beta^{-} \vee T^{-(k-2)} \alpha^{-}) \\
&= H(\alpha |T^{k} \beta^{-} \vee \alpha^{-}) + H(\alpha^{k-1} |T\beta^{-} \vee T^{-(k-2)} \alpha^{-}) \\
&= H(\alpha |T^{k} \beta^{-} \vee \alpha^{-}) + H(\alpha |T^{k-1} \beta^{-} \vee \alpha^{-}) \\
&\quad + H(\alpha^{k-2} |T\beta^{-} \vee T^{-(k-3)} \alpha^{-}) \\
&= \sum_{i=0}^{k-1}H(\alpha |T^{i+1} \beta^{-} \vee \alpha^{-}).
\end{align*}

\begin{align*}
\frac{1}{n} H(\alpha^n |T\beta^{-} \vee T^{-(n-1)} \alpha^{-}) &= \frac{1}{n} \sum_{k=0}^{n-1}H(\alpha |\alpha^{-} \vee T^{k+1} \beta^{-}) \\
&= \frac{1}{n} \sum_{k=0}^{n-1}H(\alpha |\alpha^{-} \vee T^{k+1} (\bigvee_{i=1}^{\infty} T^{-i} \beta)) \\
&= \frac{1}{n} \sum_{k=0}^{n-1}H(\alpha |\alpha^{-} \vee (\bigvee_{i=-k}^{\infty} T^{-i} \beta)) \\
&\searrow \frac{1}{n} \sum_{k=0}^{n-1}H(\alpha |\alpha^{-} \vee \beta_T) = H(\alpha |\alpha^{-} \vee \beta_T)
\end{align*}
And since $\beta \le \beta {\rm \vee }\alpha $, by Lemma 3.7's (3.3) we know
\[\frac{1}{n} H(\beta ^{n} |T^{-(n-1)} \beta ^{-} {\rm \vee }T^{-(n-1)} \alpha ^{-} )\to h(T,\beta )\] 
Therefore, letting $n\to \infty $, we have $h(\alpha {\rm \vee }\beta )=H(\alpha |\alpha ^{-} {\rm \vee }\beta _{T} )+h(\beta )$.

Lemma 3.8 is thus proved.                                                 

\bigskip
\textbf{Lemma 3.9.} If $h(\alpha {\rm \vee }\beta )=h(\alpha ,T)+h(\beta ,T)$, then
\begin{equation}
H(\alpha |\alpha ^{-} {\rm \vee }\beta _{T} )=H(\alpha |\alpha ^{-} ).\tag{3.5}    
\end{equation}

\textbf{Proof of Lemma 3.9:} From Lemma 3.8's (3.4), we have
\begin{align*}
H(\alpha |\alpha^{-} \vee \beta_{T}) &= h(\alpha \vee \beta, T) - h(\beta, T) \\
&= h(\alpha, T) \\
&= H(\alpha |\alpha^{-}).
\end{align*}

Lemma 3.9 is thus proved.    

\bigskip
\textbf{Lemma 3.10.} If $h(\alpha {\rm \vee }\beta )=h(\alpha ,T)+h(\beta ,T)$, then
\begin{equation}
H(T^{k} \beta _{n} |\alpha ^{-} )=H(T^{k} \beta _{n} |\bigwedge_{i=1}^{\infty } T^{-i} \alpha ^{-} )\tag{3.6}
\end{equation}

For every integer $k$ .

\textbf{Proof of Lemma 3.10:} Let \(\gamma \in \mathcal{I}\), then

\[H(\alpha |\alpha ^{-} )+H(\gamma |\alpha \vee \alpha ^{-} )=H(\alpha \vee \gamma |\alpha ^{-} )=H(\gamma |\alpha ^{-} )+H(\alpha |\gamma \vee \alpha ^{-} ).\]
Therefore,\[H(\alpha |\alpha ^{-} )=H(\alpha |\gamma \vee \alpha ^{-} )\Leftrightarrow H(\gamma |\alpha ^{-} )=H(\gamma |\alpha \vee \alpha ^{-} )\] 
Let \(\gamma =T^{i} (T^{k} \beta _{n} )\), by (3.5) from Lemma 3.9, we have \(H(\alpha |\alpha ^{-} )=H(\alpha |T^{i} (T^{k} \beta _{n} )\vee \alpha ^{-} )\), hence
\[H(T^{i} (T^{k} \beta _{n} )|\alpha ^{-} )=H(T^{i} (T^{k} \beta _{n} )|\alpha \vee \alpha ^{-} )\] 
Therefore,\[H(T^{k} \beta _{n} |T^{-i} \alpha ^{-} )=H(T^{k} \beta _{n} |T^{-i+1} \alpha ^{-} )\] 
holds for all \(i\), then by induction we can obtain\[H(T^{k} \beta _{n} |\alpha ^{-} )=H(T^{k} \beta _{n} |T^{-i} \alpha ^{-} )\] 
for all \(i\) as well. Letting \(i \to \infty \) and using the continuity of \(H\), we get\[H(T^{k} \beta _{n} |\alpha ^{-} )=H(T^{k} \beta _{n} |\bigwedge_{i=1}^{\infty } T^{-i} \alpha ^{-} ).\]
Lemma 3.10 is thus proved.

\bigskip
\textbf{Lemma 3.11.} $h(\alpha \vee \beta )=h(\alpha ,T)+h(\beta ,T)$, where $\alpha$ is a measurable partition that makes the system $(X,\mathcal{A},T,\mu)$ have completely positive entropy, then $\mathcal{A}$ and $\mathcal{B}$ are independent.

\textbf{Proof of Lemma 3.11:} If $(X,\mathcal{A},T,\mu)$ has completely positive entropy, then $\bigwedge_{i=1}^{\infty} T^{-i} \alpha^{-} =(X,\emptyset) \; (a.e.)$, thus

$$H(T^{k} \beta_{n} |\alpha^{-})=H(T^{k} \beta_{n})$$
\noindent this means $T^{k} \beta_{n}$ is independent of $\alpha^{-}$. Since\noindent $\gamma_{1}$ is independent of $\gamma_{1} \Rightarrow T\gamma_{1}$ and $T\gamma_{2}$ are independent.Therefore, by choosing appropriate $k,n$, the join $\bigvee_{i=-m}^{m} T^{i} \beta$ is independent of the join $\bigvee_{i=-m}^{m} T^{i} \alpha$, thus the algebra $\mathcal{B}_{m}$ corresponding to $\bigvee_{i=-m}^{m} T^{i} \beta$ is independent of the algebra $\mathcal{A}_{m}$ corresponding to $\bigvee_{i=-m}^{m} T^{i} \alpha$.
The finite algebra sequences $\mathcal{B}_{m}$ and $\mathcal{A}_{m}$ respectively increase to $\mathcal{B}$ and $\mathcal{A}$, such that $\bigvee_{m=1}^{\infty} \mathcal{B}_{m} =\mathcal{B}$ and $\bigvee_{m=1}^{\infty} \mathcal{A}_{m} =\mathcal{A}$, hence $\mathcal{B}$ and $\mathcal{A}$ are independent.

Lemma 3.11 is thus proved.

\bigskip
\textbf{Lemma 3.12.} Let $\mu \in M(T_{1}, \mu_{0})$, $\mathcal{B}_{1}$ and $\mathcal{B}_{c}$ are independent with respect to the measure $m\times \mu$ if and only if $\mu =\mu_{0}'$.

\textbf{Proof of Lemma 3.12:} We first prove that $\mathcal{B}_{1}$ and $\mathcal{B}_{c}$ are independent with respect to the measure $m\times \mu_{0}'$. If $F\in \mathcal{B}, E\in \mathcal{F}$, then
\begin{align*}
(m\times \mu_{0} ')(\pi_{1}^{-1}(E) \cap \pi^{-1}(F)) &= \int_{G_{1}}\int_{G_{2}}\chi_{F}(yx)\chi_{E}(y)dm(y)d\mu_{0}'(x) \\
&= \int_{G_{2}}\left\{\int_{G_{1}}\chi_{F}(yx)d\mu_{0}'(x) \right\} dm(y) \\
&= \int_{G_{1}}\chi_{F}(yx)d\mu_{0}'(x) \cdot \int_{G_{2}}\chi_{E}(y)dm(y) \\
&= \int_{G_{1}}\chi_{F}(x)d\mu_{0}'(x) \cdot \int_{G_{2}}\chi_{E}(y)dm(y) \\
&= \mu_{0}'(F)\cdot m(E) \\
&= (m\times \mu_{0} ')(\pi^{-1}F)\cdot (m\times \mu_{0} ')(\pi_{1}^{-1}F).
\end{align*}
Conversely, suppose $\mathcal{B}_{1}$ and $\mathcal{B}_{c}$ are independent with respect to the measure $m\times \mu$. If $F\in \mathcal{B}, E\in \mathcal{F}$, then
\begin{align*}
(m\times \mu )(\pi_{1}^{-1} (E) \cap \pi^{-1} (F)) &= (m\times \mu )(\pi_{1}^{-1} (E))\cdot (m\times \mu )(\pi^{-1} (F)) \\
&= m(E)\cdot m*\mu (F) \\
&= m(E)\cdot \mu_{0}'(F) \\
&= (m\times \mu_{0} ')(\pi_{1}^{-1} (E) \cap \pi^{-1} (F))
\end{align*}

Therefore, $m\times \mu$ and $m\times \mu_{0}'$ are consistent over $\mathcal{B}_{1} \vee \mathcal{B}_{c}$, thus consistent over $\mathcal{F}\times \mathcal{B}$.

\noindent Thus, if $F\in \mathcal{B},$ then $\mu (F)=(m\times \mu )(\pi_{2}^{-1} F)=(m\times \mu_{0} ')(\pi_{2}^{-1} F)=\mu_{0}'(F)$, hence $\mu =\mu_{0} '$.

Lemma 3.12 is thus proved.

\bigskip
\textbf{Lemma 3.13.${}^{\ }$${}^{[}$${}^{3}$${}^{]}$${}^{\ }$${}^{[}$${}^{6}$${}^{]}$} For every ergodic, measure-preserving transformation $T$ of the space $(X, \mathcal{B}, T, \mu)$ with finite entropy $h_{\mu}(T)$, there exists a finite generator $\alpha$ such that $h_{u}(T)=h_{\mu}(T,\alpha)$

\bigskip
\textbf{Lemma 3.14.${}^{\ }$${}^{[}$${}^{5}$${}^{]}$} If an automorphism $T$ of $G$ is ergodic, then it has completely positive entropy.

\bigskip
\textbf{Theorem 3.1.} If $G_{2}$ acts freely on $G_{1}$, and and the homeomorphism $T_{1}$ $T_{2} {\rm -}$ commutes with $G_{2}$, where $T_{2} :G_{2} \to G_{2} $ is an isomorphism, if $\mu _{0} \in M(T_{0})$, then:

\begin{enumerate}
\item[(1)]   $h(T_{1},\mu )\le h(T_{1},\mu _{0} ')\; (\mu \in M(T_{1},\mu _{0}))$;

\item[(2)]  If $\mu _{0} '$ is ergodic with respect to $T_{1}$, the Haar measure $m$ on $G_{2}$ is ergodic with respect to $T_{2}$, and $h(T_{1},\mu _{0} ')<\infty $, then
\[h(T_{1},\mu )<h(T_{1},\mu _{0} ')\; (\mu \in M(T_{1},\mu _{0}), \mu \ne \mu _{0} ').\] 
\end{enumerate}

\textbf{Proof of Theorem 3.1:} If $h(T_{2}, m) = \infty$, then $h(T, \mu_{0}') = \infty$, and the conclusion holds.

\noindent Suppose $h(T_{2}, m) < \infty$, then
\begin{align*}
h(T_{1}, \mu) &= h(T_{2} \times T_{1}, \mathcal{B}\times \mathcal{B}, m\times \mu) - h(T_{2}, m) \notag \\
&= h(T_{2} \times T_{1}, \mathcal{B}_{1} \times \mathcal{B}_{c}, m\times \mu) - h(T_{2}, m) \notag \\
&\le h(T_{2} \times T_{1}, \mathcal{B}_{1}, m\times \mu) + h(T_{2} \times T_{1}, \mathcal{B}_{c}, m\times \mu) - h(T_{2}, m) \notag \\
&= h(T_{1}, m*\mu) + h(T_{2}, m) - h(T_{2}, m) \notag \\
&= h(T_{1}, \mu_{0}').
\end{align*}

If $\mu \in M(T_{1}, \mu_{0})$ is an ergodic measure,

\noindent For contradiction, if $h(T_{1}, \mu) = h(T_{1}, \mu_{0}'){\rm \; }$, then
\begin{align*}
h(T_{2}, \mathcal{B}_{1}, m) + h(T_{1}, \mathcal{B}_{c}, \mu) &= h(T_{2} \times T, \mathcal{B}\times \mathcal{B}, m\times \mu) \notag \\
&= h(T_{2} \times T_{1}, \mathcal{B}_{1} \vee \mathcal{B}_{c}, m\times \mu) \notag \\
&\le h(T_{2} \times T_{1}, \mathcal{B}_{1}, m\times \mu) + h(T_{2} \times T_{1}, \mathcal{B}_{c}, m\times \mu) \notag \\
&= h(T_{2}, \mathcal{F}, m) + h(T_{1}, \mathcal{B}_{c}, m*\mu) \notag \\
&= h(T_{2}, \mathcal{F}, m) + h(T_{1}, \mathcal{B}, \mu_{0}').
\end{align*}

Since both sides of the equation are equal, then
\[h(T_{2} \times T_{1}, {\rm {\mathcal B}}_{1} {\rm \vee {\mathcal B}}_{c}, m\times \mu) = h(T_{2} \times T_{1}, {\rm {\mathcal B}}_{1}, m\times \mu) + h(T_{2} \times T_{1}, {\rm {\mathcal B}}_{c}, m\times \mu)\]
Since the systems $(T_{2} \times T_{1}, \mathcal{B}_{1}, m\times \mu)$ and $(T_{2} \times T_{1}, \mathcal{B}_{c}, m\times \mu)$ are respectively isomorphic to $(T_{2}, \mathcal{F}, m)$ and $(T_{1}, \mathcal{B}, \mu_{0}')$, they are ergodic and have finite entropy. By Lemma 3.13, $(T_{2} \times T_{1}, \mathcal{B}_{1}, m\times \mu)$ and $(T_{2} \times T_{1}, \mathcal{B}_{c}, m\times \mu)$ have finite entropy generators $\alpha$ and $\beta$, thus $\mathcal{B}_{1} = \mathcal{A}, \mathcal{B}_{c} = \mathcal{B}\; a.e.$, and by Lemma 3.14, $(T_{2} \times T_{1}, \mathcal{B}_{1}, m\times \mu)$ has completely positive entropy. Combined with Lemma 3.11, it follows that $\mathcal{B}_{1}$ and $\mathcal{B}_{c}$ are independent with respect to $m\times \mu$. By Lemma 3.12, $\mu = \mu_{0}'$, contradicting the assumption.

Therefore, if $\mu \ne \mu_{0}', \mu \in M(T_{1}, \mu_{0})$ and $\mu$ is ergodic, then $h(T_{1}, \mu) < h(T_{1}, \mu_{0}')$. By Corollary 3.1, if $\mu \in M(T_{1}, \mu_{0}), \mu \ne \mu_{0}'$, then $h(T_{1}, \mu) < h(T_{1}, \mu_{0}')$.

Theorem 3.1 is thus proved.

\bigskip
\textbf{Corollary 3.2.}  If $G_{2}$ acts freely on $G_{1}$, and and the homeomorphism $T_{1}$ $T_{2} {\rm -}$ commutes with $G_{2}$, where $T_{2} :G_{2} \to G_{2} $ is an isomorphism, assuming $T_{0}$ has a unique measure of maximal entropy $\mu_{0}$,and if $\mu_{0}'$ is ergodic with respect to $T$, the Haar measure $m$ on $G_{2}$ is ergodic with respect to $T_{2}$, and $h(T_{1}, \mu_{0}')<\infty$, then $T_{1}$ has a unique measure of maximal entropy $\mu_{0}'$.

\textbf{Proof of Corollary 3.2:} Since $h(T_{1}, \mu_{0}') < \infty$, we know $h(T_{2}, m) < \infty$. Let $\nu_{0} \in M(T_{0})$, for $\nu \in M(T_{1}, \nu_{0})$, by the first part of Theorem 3.1 and Lemma 3.3, we have
\[h(T_{1}, \nu) \le h(T_{1}, \nu_{0}') = h(T_{0}, \nu_{0}) + h(T_{2}, m)\]

Let $h(T_{1} ,\nu )=h(T_{1} ,\mu _{0} ')$ and $\nu \in M(T_{1} ,\nu _{0} )$, since $h(T_{1} ,\mu _{0} ')$ is the maximum entropy, thus by
$h(T_{1} ,\nu )\le h(T_{1} ,\nu _{0} ')\le h(T_{1} ,\mu _{0} ')$
and $h(T_{1} ,\nu )=h(T_{1} ,\mu _{0} ')$,we have $h(T_{1} ,\nu )=h(T_{1} ,\nu _{0} ')=h(T_{0} ,\nu _{0} )+h(T_{2} ,m)$,and since $h(T_{1} ,\mu _{0} ')=h(T_{0} ,\mu _{0} )+h(T_{2} ,m)$, therefore
\[\begin{array}{l} {h(T_{0} ,\mu _{0} )+h(T_{2} ,m)=h(T_{1} ,\mu _{0} ')=h(T_{1} ,\nu )=h(T_{0} ,\nu _{0} )+h(T_{2} ,m)} \\ {\Rightarrow h(T_{0} ,\mu _{0} )=h(T_{0} ,\nu _{0} )} \end{array}\] 

Since $T_{0} $ has a unique maximum entropy measure $\mu _{0} $, thus $\nu _{0} =\mu _{0} $. Therefore, $\nu \in M(T_{1} ,\mu _{0} )$, according to the second part of the theorem 3.1, if $h(T_{1} ,\nu )=h(T_{1} ,\mu _{0} ')$, then $\nu =\mu _{0} '$.

Corollary 3.1 is thus proved.

\bigskip
\textbf{Theorem 3.2.} For every affine transformation $T=aA$ of a compact metrizable group $G$, we have

\noindent $h(T,\mu) \le h(T,m) \; (\mu \in M(T))$. If $m$ is ergodic with respect to $T$, and $h(T,m)<\infty$, $h(T,\mu)<h(T,m) \; (\mu \in M(T), \mu \ne m)$

\textbf{Proof of Theorem 3.2:} $G$ acts on itself by left multiplication, and$T(yx)=aA(yx)=aA(y)A(x)=aA(y)a^{-1} T(x)$. Let $B(y)=aA(y)a^{-1} $, then $B$ is also an automorphism of $G$, and $T(yx)=B(y)T(x)$, which follows from Theorem 3.1.

Theorem 3.2 is thus proved.

\bigskip
\textbf{Lemma 3.15.${}^{[8]}$ } Properties of the natural extension:
\begin{enumerate}
\item[(1)] The natural extension of a measure-preserving continuous mapping $T$ is a measure-preserving automorphism.
\item[(2)]  The entropy of $T$ is equal to the entropy of its natural extension, that is, $h(T,\mu )=h(S,\underline{\mu })$.
\end{enumerate}
\bigskip

For non-invertible maps, we have

\textbf{Theorem 3.3.} Let $G_{2}$ be a compact metrizable group acting freely on $G_{1}$, and the continuous map $T_{1}:G_{1} \to G_{1}$ $T_{2}-$commutes with $G_{2}$  where $T_{2}: G_{2} \to G_{2}$ is a homomorphism. Assuming $\mu_{0} \in M(T_{0})$, then:

\begin{enumerate}
\item[(1)] $h(T_{1}, \mu) \le h(T_{1}, \mu_{0}') \; (\mu \in M(T_{1}, \mu_{0}))$;

\item[(2)] If $\mu_{0}'$ is ergodic with respect to $T_{1}$, the Haar measure $m$ on $G_{2}$ is ergodic with respect to $T_{2}$, and $h(T_{1}, \mu_{0}')<\infty$, then
\[h(T_{1}, \mu)<h(T_{1}, \mu_{0}') \; (\mu \in M(T_{1}, \mu_{0}), \mu \ne \mu_{0}').\]
\end{enumerate}

\textbf{Proof of Theorem 3.3:} Let $Y$ be the compact metrizable group defined by
\[Y=\{ (x_{1} ,x_{2} ,x_{3} ,\ldots )|Tx_{i+1} =x_{i} ,i=1,2,\ldots \} \]

and let $S$ be the homeomorphism of $Y$ defined by $S(x_{1} ,x_{2} ,\ldots )=(Tx_{1} ,x_{1} ,x_{2} ,\ldots )$. Let $H$ be the compact metrizable group defined by
\[\{ (y_{1} ,y_{2} ,\ldots )|T_{2} y_{i+1} =y_{i} ,i=1,2,\ldots \} \] 

and let $B$ be the automorphism on $H$ defined by $B(y_{1} ,y_{2} ,\ldots )=(T_{2} y_{1} ,y_{1} ,y_{2} ,\ldots )$.

$H$ acts freely on $Y$ by $((y_{1} ,y_{2} ,\ldots ),(x_{1} ,x_{2} ,\ldots ))\to ((y_{1} x_{1} ,y_{2} x_{2} ,\ldots )$, and $S$ commutes with $H$ under $B$-action. Furthermore, the orbit space $Y_{0} =Y/H$ corresponds one-to-one with
\[
\left\{ (z_{1}, z_{2}, \ldots) \mid z_{i} \in G_{0},\ T_{0}z_{i+1} = z_{i},\ i \geq 1 \right\}
\]

and the induced homeomorphism $S_{0} :Y_{0} \to Y_{0} $ is given by $S_{0} (z_{1} ,z_{2} ,\ldots )=(T_{0} z_{1} ,z_{1} ,z_{2} ,\ldots )$. If $\nu $ is a measure on $G_{1} $ (or $G_{0} $), let $\underline{\nu }$ represent the corresponding inverse limit measure on $Y$ (or $Y_{0} $). By Lemma 3.15, then$h(T_{1} ,\nu)=h(S,\underline{\nu })$.Furthermore, the Haar extension of $\underline{\mu _{0} }$ is the inverse limit measure corresponding to $\mu _{0} '$, denoted by $\underline{\mu _{0} '}$. Let $\mu \in M(T_{0} ,\mu _{0} )$.Then $\underline{\mu }\in M(S_{0} ,\underline{\mu _{0} })$, and by the first part of Theorem 3.1, we have $h(S,\underline{\mu })\le h(S,\underline{\mu _{0} '})$. Therefore$h(T_{1} ,\mu )\le h(T_{1} ,\mu _{0} ').$

Now assume that $\mu _{0} '$ is an ergodic measure with respect to $T_{1} $ and $h(T_{1} ,\mu _{0} ')<\infty $. Then $\underline{\mu _{0} }$ is ergodic with respect to $S$, and $h(S,\mu _{0} ')=h(T_{1} ,\mu _{0} ')<\infty $. If $\mu \in M(T_{1} ,\mu _{0} ){\rm \; }(\mu \ne \mu _{0} ')$ then $\mu \in M(S,\mu _{0} ){\rm \; }(\mu \ne \mu _{0} ')$, and by Theorem 3.1, we have $h(S,\mu )<h(S,\mu _{0} ')$. Therefore
$h(T_{1} ,\mu )=h(S,\underline{\mu })<h(S,\underline{\mu _{0} '})=h(T_{1} ,\mu _{0} ').$

Theorem 3.3 is thus proved.

\bigskip
\textbf{Corollary 3.3:} If $G_{2}$ acts freely on $G_{1}$, and $T_{1}:G_{1} \to G_{1}$ $T_{2}-$commutes with $G_{2}$. Assume $T_{0}$ has a unique measure $\mu_{0}$ with maximal entropy . If $\mu_{0}'$ is an ergodic measure with respect to $T_{1}$, and $h(T_{1}, \mu_{0}')<\infty$, then $T_{1}$ has a unique measure $\mu_{0}'$ with maximal entropy .

\textbf{Proof of Corollary 3.3:} The proof of this corollary is similar to the proof of Corollary 3.2, since $h(T_{1} ,\mu _{0} ')=h(T_{0} ,\mu _{0} )+h(T_{2} ,m)$ also holds for the homomorphism $T_{2} $.

Corollary 3.3 is thus proved.

\bigskip
\textbf{Corollary 3.4:} Let $G$ be a compact metrizable group, $T:G\to G$ a  homomorphism, $m$ is the Haar measure on $G$, then:

\begin{enumerate}
\item[(1)] $h(T,\mu) \le h(T,m) \; (\mu \in M(T))$;

\item[(2)] If $m$ is ergodic with respect to $T$, and $h(T,m)<\infty$, then $h(T,\mu)<h(T,m) \; (\mu \in M(T), \mu \ne m).$
\end{enumerate}

\vspace{1cm}
\begin{center}
{\bf\large 4. Ergodicity of $\mu\ast\nu$}
\end{center}

\textbf{Theorem 4.1.} Let $T$ be a surjective homomorphism on $G$. If $(G,T,{\mathcal F},\mu)$ and $(G,T,{\mathcal F},\nu)$ are disjoint ergodic dynamical systems, then $\mu * \nu$ is ergodic.

\textbf{Proof of Theorem 4.1:} By contradiction, suppose $\mu *\nu $ is not ergodic, then there exists a measurable set $E\in {\rm {\mathcal B}}$ satisfying $T^{-1} E=E$ such that $0<\mu *\nu (E)<1$. Since $T$ is surjective, we have $E=T(T^{-1} E)=TE$, define a measure on $G\times G$ as $\lambda $:
\[\lambda (F)=\frac{1}{\mu *\nu (E)} \mu \times \nu (F\cap \pi ^{-1} (E)){\rm \; (}F\in {\rm {\mathcal B}}\times {\rm {\mathcal B}}).\] 
Let $A\in {\rm {\mathcal B}}$, define a function on $G\times G$ as $f(x,y)=\chi _{A} (x)$, where $\chi _{A} (x)$ represents the characteristic function of $A$. Thus $d\lambda =\frac{\chi _{\pi ^{-1} E} d(\mu \times \nu )}{\mu *\nu (E)} $, therefore
\begin{align*}
\int _{G\times G}f(x,y)d\lambda &= \int _{G\times G}f(x,y)\frac{\chi _{\pi ^{-1} E} (x,y)}{\mu *\nu (E)}  d(\mu \times \nu ) \\
&= \frac{1}{\mu *\nu (E)} \int _{G}\int _{G}\chi _{A}  (x)\chi _{\pi ^{-1} E} (x,y)d\nu (y)d\mu (x)  \\
&= \frac{1}{\mu *\nu (E)} \int _{G}\chi _{A} (x)\int _{G}\chi _{\pi ^{-1} E} (x,y)d\nu (y)d\mu (x).
\end{align*}
Considering the inner integral

\noindent $\int _{G}\chi _{\pi ^{-1} E} (x,y)d\nu (y) =\nu (\{ y:xy\in E\} )$.

Because
\[\begin{array}{l} {\{ y:xy\in E\} =x^{-1} \cdot E=T^{-1} [T(x^{-1} \cdot E)]=T^{-1} [T(x^{-1} )\cdot TE]} \\ {=T^{-1} [T(x^{-1} )\cdot E]=T^{-1} [(Tx)^{-1} \cdot E]=T^{-1} \{ y:(Tx)y\in E\} } \end{array}\] 

Since $\nu $ is invariant under $T$, we have
\[\nu (\{ y:xy\in E\} =\nu (T^{-1} \{ y:(Tx)y\in E\} )=\nu (\{ y:(Tx)y\in E\} )\] 

Therefore
\[\int _{G}\chi _{\pi ^{-1} E} (x,y)d\nu (y)= \int _{G}\chi _{\pi ^{-1} E} (Tx,y)d\nu (y) \] 

Since $T$ is surjective, for every $x'\in G, \exists x\in G$, so $\int _{G}\chi _{\pi ^{-1} E} (x,y)d\nu (y) $ is a function invariant in $x$, therefore $\int _{G}\chi _{\pi ^{-1} E} (x,y)d\nu (y)= K{\rm \; \; a.e}$, $K$ is a constant. Therefore
\[\int _{G\times G}f(x,y)d\lambda = \frac{K}{\mu *\nu (E)} \int _{G}\chi _{A} (x)d\mu (x)= \frac{K\mu (A)}{\mu *\nu (E)} .\] 

Letting $A=G$, then $K=\mu *\nu (E)$. Hence
\[\int _{G\times G}f(x,y)d\lambda = \mu (A)\] 

This shows that, if $f(x,y)$ depends only on $x$, then
\[\int _{G\times G}f(x,y)d\lambda = \int _{G\times G}f(x,y)d\mu \circ \pi ^{-1}  \] 

The above process is also applicable to $f(x,y),y$ and $\nu $. Combining this fact, and that $\lambda \ne \mu \times \nu $, means that $\pi _{1} {}^{-1} ({\rm {\mathcal B}})$ and $\pi _{2} {}^{-1} ({\rm {\mathcal B}})$ are not independent. In fact, if $\pi _{1} {}^{-1} ({\rm {\mathcal B}})$ and $\pi _{2} {}^{-1} ({\rm {\mathcal B}})$ were independent, then for every $F\in {\rm {\mathcal B}}$, we have
\begin{align*}
\lambda (\pi _{1} {}^{-1} (F)\cap \pi _{2} {}^{-1} (F)) & = \lambda (\pi _{1} {}^{-1} (F))\cdot \lambda (\pi _{2} {}^{-1} (F)) \\
& = \lambda (\{ (x,y):x\in F\} )\cdot \lambda (\{ (x,y):y\in F\} ) \\
& = \int _{G\times G}\chi _{F} (x)d\lambda \cdot  \int _{G\times G}\chi _{F} (y)d\lambda \\
& = \mu (F)\cdot \nu (F) \\
& = \mu \times \nu (\{ (x,y):x\in F\} )\cdot \mu \times \nu (\{ (x,y):y\in F\} ) \\
& = \mu \times \nu (\pi _{1} {}^{-1} F)\cdot \mu \times \nu (\pi _{2} {}^{-1} F) \\
& = \mu \times \nu (\pi _{1} {}^{-1} F\cap \pi _{2} {}^{-1} F).
\end{align*}
This contradicts $\lambda \ne \mu \times \nu$.

By the definition of disjointness, given $(G,T,{\rm {\mathcal B}},\mu )=X,(G,T,{\rm {\mathcal B}},\nu )=Y,\pi _{1} =\alpha ,\pi _{2} =\beta $, and $(G\times G,T\times T,{\rm {\mathcal B}}\times {\rm {\mathcal B}},\mu \times \nu )=Z$, we know that $(G,T,{\rm {\mathcal B}},\mu )$ and $(G,T,{\rm {\mathcal B}},\nu )$ are not disjoint, which contradicts the given condition. Therefore, $\mu *\nu$ is exhaustive.

Theorem 4.1 is thus proved.

\vspace{1cm}
\begin{center}
{\bf\large Conclusion}
\end{center}

This paper discusses the maximum entropy theorem, invariant measure convolutions' ergodicity, and extends the application of the maximum entropy theorem and measure convolution ergodicity from isomorphism to surjective homomorphism. This expansion provides new tools and perspectives for understanding complex behaviors in dynamical systems.

Theorem 3.1 examines the relationship between a measure's entropy in a dynamical system and its maximum entropy under certain conditions. It assumes a scenario where one group freely acts on another and connects the actions of these two groups through a specific function. The theorem states that the entropy of any measure does not exceed that of a certain maximal entropy measure. Furthermore, if additional conditions of ergodicity and finiteness are met, then the entropy of all other measures is strictly less than the maximal entropy, except for the maximal entropy measure itself, indicating the existence and uniqueness of the system's maximal entropy measure.

Theorem 3.3 expands the conclusions of Theorem 3.1 under homomorphisms, offering new insights into the ergodicity of dynamical systems, especially when considering the convolution of measures from two disjoint ergodic systems. It shows that under specific conditions, the convolution of two ergodic systems is ergodic, providing an essential theoretical foundation for understanding the ergodic behavior of composite systems.

Overall, these findings not only deepen our understanding of the maximum entropy theorem and ergodicity theory but also significantly broaden their applicability by extending these concepts from isomorphic to surjective homomorphic dynamical systems. This has important theoretical and practical implications for studying invariant measures and ergodic properties in dynamical systems and understanding how these systems interact within a broader mathematical framework.
\vspace{1cm}
\begin{center}
{\bf\large References}
\end{center}

\let\oldsection\section
\renewcommand{\section}[2]{}%

\let\section\oldsection

\end{document}